\documentclass[12pt,a4paper]{amsart}
\usepackage{amssymb,amsfonts, amsmath, mathrsfs}

\title[]{Lower bounds for the first Laplacian eigenvalue of geodesic balls of spherically
symmetric manifolds}
\author{Cleon S. Barroso and G. Pacelli  Bessa}
\thanks{The first author is grateful for the financial support by CAPES - PRODOC}
\address{Departamento de Matem\'atica, Universidade Federal do Cear\'a,
Campus do Pici, Bl. 914, 60455-760, Fortaleza, CE, Brazil.}
\email{cleonbar@mat.ufc.br, bessa@mat.ufc.br}


\newtheorem{theorem}{Theorem}[section]
\newtheorem{corollary}{Corollary}[section]
\newtheorem{lemma}{Lemma}[section]
\newtheorem{definition}{Definition}[section]

\theoremstyle{remark}
\newtheorem{remark}[theorem]{Remark}

\allowdisplaybreaks
\begin{document}
\begin{abstract}We  obtain  lower bounds
for the first Laplacian eigenvalues of geodesic balls of spherically
symmetric manifolds. These lower bounds are only $C^{0}$ dependent
on the metric coefficients.

\end{abstract}
\maketitle
\section{Introduction}Let
 $B(r)$ be a geodesic ball of radius
    $r$ in  the
$n$-dimensional  sphere $\mathbb{S}^{\,n}(1)$  of sectional
curvature $+1$.  Although the sphere is a well studied manifold, the
values of the  first Laplacian eigenvalue $\lambda_{1}(r)$  on
$B(r)$,  (Dirichlet boundary data if $r<\pi$) are pretty much
unknown, exceptions are $\lambda_{1}(\pi/2)=n$ and
$\lambda_{1}(\pi)=0$. Among the  various types of bounds for
$\lambda_{1}(r)$, see \cite{Barbosa-do-Carmo}, \cite{pinsky},
\cite{Sato} in dimension two, see \cite{friedland-hayman} in
dimension three, we would like to emphasize the following bounds due
to
 Betz, Camera
and Gzyl they  obtained in \cite{betz-camera-gzyl}.
\begin{equation}\label{eqBCG}\left(\frac{c_{n}}{r}\right)^{2}>\lambda_{1}(r)\geq
\frac{1}{\int_{0}^{r}[\frac{1}{\sin^{n-1}(\sigma)}\cdot
\int_{0}^{\sigma}\sin^{n-1} (s)ds]\,d\sigma}, \end{equation}Where
$c_{n}$ is the first zero of the $J_{(n-2)/2}$ Bessel function. The
upper bound is just Cheng's eigenvalue comparison theorem
 \cite{cheng2} and it is due to the fact that
the Ricci curvature of the sphere is positive (need only to be
non-negative). The interesting part is the lower bound that they
obtained with probabilistic method. Denoting by $V(r)$ the
$n$-volume of the geodesic ball $B(r)$ and by $S(r)$ the
$(n-1)$-volume of the boundary $\partial B(r)$ we can rewrite
Betz-Camera-Gzyl lower bound as
\begin{equation}\label{eqBCG2}\lambda_{1}(r)\geq
\frac{1}{\displaystyle\int_{0}^{r}\frac{V(\sigma)}{S(\sigma)}
\,d\sigma}.
\end{equation}In this note, using a fixed point theorem approach, we extend Betz-Camera-Gzyl's lower bound
to $\lambda_{1}(r)$ of geodesic balls
 $B(r)$  of complete spherically symmetric manifolds.

 \noindent  A
 spherically symmetric manifold  is a quotient space $M=([0,R)\times
  \mathbb{S}^{n-1})/\backsim$, with $R\in (0,\infty]$, where $$(t,\theta)\backsim (s,\alpha)
\Leftrightarrow \left\{\begin{array}{lll} t=s \,\,\,\,\, and &&
\theta =\alpha
  \\
or\,\,& & \\s=t=0.&&\end{array}\right.$$ endowed with a Riemannian
metric of this form
 $dt^{2}+f^{2}(t)d\theta^{2}$, $f(0)=0$, $f'(0)=1$, $f(t)>0$ for
all $t\in (0,R]$.    The class of spherically symmetric manifolds
includes the canonical space forms $\mathbb{R}^{n}$,
$\mathbb{S}^{n}(1)$ and $\mathbb{H}^{n}(-1)$. A spherically
symmetric manifold has a pole (at $p=\{0\}\times \mathbb{S}^{n-1}$)
if and only if $R=\infty$.

\begin{theorem}\label{thmBarroso-Bessa}Let $M=[0,R)\times
  \mathbb{S}^{n-1}$ be a spherically symmetric manifold  with Riemannian metric
$dt^{2}+f^{2}(t)d\theta^{2}$, $f(0)=0$, $f'(0)=1$, $f(t)>0$ for all
$t\in (0,R]$ and $B(r)\subset M$  a geodesic ball of radius $r$.
Then \begin{equation} \lambda_{1}(r)\geq
\frac{1}{\displaystyle\int_{0}^{r}\frac{V(\sigma)}{S(\sigma)}
\,d\sigma}.\label{eqBarroso-Bessa}
\end{equation}
\end{theorem}

\begin{definition}Let $M$ be a spherically symmetric manifold with a
pole. The fundamental tone $\lambda^{\ast}(M)$ is defined by
\begin{equation}\lambda^{\ast}(M)=\lim_{r\to \infty}\lambda_{1}(r)
\end{equation}
\end{definition}

\begin{corollary}\label{Cor-Barroso-Bessa}Let $M=[0,\infty)\times
  \mathbb{S}^{n-1}$ be a spherically symmetric manifold with a pole.
  Then \begin{equation}\lambda^{\ast}(M)\geq
\frac{1}{\displaystyle\int_{0}^{\infty}\frac{V(\sigma)}{S(\sigma)}
\,d\sigma}.\label{eqBarroso-Bessa2}\end{equation}
\end{corollary}
This corollary is closely related to  certain property of the
Brownian motions on M.  Denote by $p(t,x,y)\in
C^{\infty}((0,\infty)\times M\times M)$  the heat kernel of $M$ and
let $X_{t}$ be a Brownian motion on $M$ and denote by
$\mathbb{P}_{x}$ the corresponding measure in the space of paths
emanating from a point $x$. See more details in \cite{grigor'yan}.
\begin{definition}A Brownian motion $X_{t}$ on a complete manifold
$M$ is recurrent if for any $x\in M$ and any non-empty open set
$\Omega \subset M$ \begin{equation}\mathbb{P}_{x}\left(\{There
\;\,is\;\,a\;\, sequence\; \,t_{k}\to \infty\,such\;\, that\;\,
X_{t_{k}}\in \Omega\}\right)=1.
\end{equation}Otherwise is transient.
\end{definition}
\begin{definition}A Brownian motion $X_{t}$ on a complete manifold
$M$ is stochastically complete  if for all $x\in M$ and $t>0$.

\begin{equation}\int_{M}p(t,x,y)d\mu (y)=1
\end{equation} Otherwise $X_{t}$ is incomplete.
\end{definition}

We say that a complete  manifold $M$ is recurrent, transient,
stochastically complete, incomplete if the Brownian motion has this
property. The following test is well known, see \cite{grigor'yan}
and references there in.

\vspace{3mm}
 \noindent {\bf Test for Stochastically Completeness:} {\em Let $M$ a
spherically symmetric manifold with a pole. Then $M$ is
stochastically complete if and only if
$$\int_{0}^{\infty}\frac{V(r)}{S(r)}dr=\infty.
$$}

\begin{remark}\begin{itemize}\item[]
\item[i.] Let $M$ be a complete Riemannian manifold. If
$\lambda^{\ast}(M)>0$ then  $M$ is transient.
\item[ii.] There are examples of complete, stochastically incomplete
(therefore transient) Riemannian manifolds $M$ with
$\lambda^{\ast}(M)=0$, see \cite{pinchover}.
\end{itemize}
\end{remark}The following corollary follows from the test for
stochastically completeness and Corollary (\ref{Cor-Barroso-Bessa}).
\begin{corollary}Let $M$ be a
spherically symmetric manifold with a pole. If $M$ is stochastically
incomplete then $\lambda^{\ast}(M)>0$. If $\lambda^{\ast}(M)=0$ then
$M$ is stochastically complete.
\end{corollary}
\section{Proof of the results}
Consider the space $X$ of all continuous functions on $[0,r]$ with
the usual topology defined by the norm $\|u\|=\sup\sb{0\leq t\leq
r}|u(t)|$. For  $a\in \mathbb{R} $ and $\Theta>0$  let $T=T_{a,
\Theta}$ be the operator in $X$ defined by
$$
T\,u(t)=\Theta-\int\sb 0\sp t\int\sb
0\sp\sigma\Big(\frac{f^{n-1}(s)}{f^{n-1}(\sigma)}\Big)[a+\lambda_{1}(r)]u(s)\,ds\,d\sigma,\quad
0\leq t\leq r
$$ Let $B(r)\subset M$ be a geodesic ball of radius $r<R$ in a spherically symmetric manifold
 $M=[0,R)\times \mathbb{S}^{n-1}$
with metric $dt^{2}+f^{2}(t)d\theta^{2}$.  The Laplacian operator
$\triangle_{M}$ at a point $(t,\theta)$ is given by
$$
\triangle_{M} =\frac{\partial^{2}}{\partial_
{t}}+(n-1)\frac{f'(t)}{f(t)}\frac{\partial}{\partial_{t}}
+\frac{1}{f^{2}(t)}\triangle_{\mathbb{S}^{n-1}}$$

\noindent Given $u\in X$, we can extend (radially) $u$ and $Tu$ to
continuous functions $\tilde{u}$ and $\tilde{T}u$ on $B(r)$
respectively by $\tilde{u}(t,\theta)=u(t)$ and
 $\tilde{T}u(t,\theta)=Tu(t)$,
for all $\theta \in \mathbb{S}^{n-1}$, $t\in [0, r)$. A straight
forward computation shows that \begin{equation}\triangle \tilde{T}u
(t,\theta)
+(a+\lambda_{1}(r))\,\tilde{u}(t,\theta)=0\label{eqBarroso-Bessa3}\end{equation}
for all $t\in[0,r]$ and all $\theta\in \mathbb{S}^{n-1}.$

\vspace{2mm} \noindent Let
$C(r)=\displaystyle\int_{0}^{r}\left[\frac{1}{f^{n-1}(\sigma)}\int_{0}^{\sigma}f^{n-1}(s)ds\right]d\sigma
=\int_{0}^{r}\frac{V(\sigma)}{S(\sigma)}d\sigma$. Suppose that
$\lambda_{1}(r)< C(r)^{-1}$ and choose $a>0$  such that
$\lambda_{1}(r)+a< C(r)^{-1}$. We will show  that the operator
$T_{a, \Theta}$ has a fixed point $u_{a,\Theta}$ in the closed
convex subset $F=\{u\in X\colon 0\leq u\leq \Theta\}$ of $X$. If
$Tu_{a,\Theta}=u_{a,\Theta}$ then the radial extensions
$\tilde{u}_{a,\Theta}$ and $\tilde{T}u_{a,\Theta}$ satisfies by
(\ref{eqBarroso-Bessa3}) the following identity.
\begin{equation}\label{eqBarroso-Bessa4}\triangle \tilde{u}_{a,\Theta} (t,\theta)
+(a+\lambda_{1}(r))\,u_{a,\Theta}(t,\theta)=0\end{equation} for all
$t\in[0,r]$ and all $\theta\in \mathbb{S}^{n-1}$. But this
contradicts the following well known lemma.
\begin{lemma}There is no non-trivial smooth  solution to the
problem $$\left\{\begin{array}{rlll} \triangle u +
(a+\lambda_{1}(r))u &=&0 & in \,\,\, B(r)\\
u&\geq &0 & in \,\,\,\overline{B(r)},
\end{array}\right.$$if $a>0$.
\end{lemma}\noindent Thus we have that $\lambda_{1}(r)\geq
C(r)^{-1}$, proving (\ref{eqBarroso-Bessa}).

\vspace{.2cm} To finish the proof of Theorem
(\ref{thmBarroso-Bessa}) we need to show that $T_{a,\Theta}:F\to F$
has a fixed point.
 In order to get a fixed point for $T_{a,\Theta}$, we are going to use the following
well known Schauder-Tychonoff fixed point theorem.

\begin{theorem}\label{trm:fp} Let $F$ be a nonempty closed convex subset
of a separated locally convex topological vector space $X$. Suppose
that $T\colon F\to F$ is a continuous map such that $T(F)$ is
relatively compact. Then $T$ has a fixed point.
\end{theorem}

We are going to show that $T_{a,\Theta}$ satisfies the hypotheses of
Theorem (\ref{trm:fp}) if $\lambda_{1}(r)+a< C(r)^{-1}$. We start we
few lemmas.

\begin{lemma}\label{lem:2} Let
$F$ be the set
$$
F=\{u\in X\colon 0\leq u(r)\leq \Theta\}
$$
Then $T$ maps $F$ into itself.
\end{lemma}
\begin{proof} Let $u\in F$ be arbitrary.
Clearly, $Tu$ is continuous. Since $(a+\lambda_{1})u\geq 0$, we have
that $\int\sb 0\sp t\int\sb
0\sp\sigma\Big(\frac{f^{n-1}(s)}{f^{n-1}(\sigma)}\Big)[a+\lambda_{1}(r)]u(s)\,ds\,d\sigma
\geq 0$ thus
 $(Tu)(t)\leq \Theta$, for all
$0\leq t <r$. On the other hand, since
$(a+\lambda_{1}(r))<C(r)^{-1}$ and $0\leq u(t)\leq \Theta$, we have
that,
\begin{eqnarray}\label{eqn:positivity}
(Tu)(t)&=&  \Theta - \int\sb 0\sp t\int\sb
0\sp\sigma\Big(\frac{f^{n-1}(s)}{f^{n-1}(\sigma)}\Big)[a+\lambda_{1}(r)]u(s)\,ds\,d\sigma
\nonumber
\\
&\geq& \Theta - \int\sb 0\sp r\int\sb
0\sp\sigma\Big(\frac{f^{n-1}(s)}{f^{n-1}(\sigma)}\Big)[a+\lambda_{1}(r)]u(s)\,ds\,d\sigma
\nonumber\\
&\geq & \Theta - \int\sb 0\sp r\int\sb
0\sp\sigma\Big(\frac{f^{n-1}(s)}{f^{n-1}(\sigma)}\Big)C^{-1}(r)\,\Theta\,ds\,d\sigma\nonumber
\\
&=& 0\nonumber
\end{eqnarray}
for all $0\leq t<r$. This proves that $T(F)\subset F$.
\end{proof}
\begin{lemma}The map $T=T_{a,\Theta}\colon F \to F$ is
continuous and $T(F)$ is relatively compact.
\end{lemma}
\begin{proof}Note that $F$ is closed and convex. Let
$\{u\sb m\}\subset F$ be a sequence such that $u\sb m\to u$, for
some $u\in F$, (recall that $\Vert u \Vert = \sup_{0\leq s\leq r}
\vert u(s)\vert $).  Thus, we have
\begin{eqnarray*}
|Tu\sb n(t)-Tu(t)|\leq \Vert u_{n}-u\Vert \, [a+\lambda_{1}(r)]
\int\sb 0\sp t\int\sb
0\sp\sigma\Big(\frac{f^{n-1}(s)}{f^{n-1}(\sigma)}\Big)\,ds\,d\sigma.
\end{eqnarray*}We can conclude
that $Tu\sb m$ converges uniformly to $Tu$. Moreover, $$ \vert
Tu'(t)\vert \leq \frac{\Theta\, C^{-1}(r)}{f^{n-1}(t)}
\int_{0}^{t}f^{n-1}(s)ds=h(t)
$$Observe that $h(t)$ is a continuous function on $[0,r]$ thus $ \vert
Tu'(t)\vert \leq \sup_{[0,r]}h(t)$ which implies that each $T(F)$ is
equicontinuous.   Since $T(F)$ is uniformly bounded, the
Ascoli-Arzela theorem implies that $T(F)$ is relatively compact.
\end{proof}

\end{document}